\documentclass[a4paper,12pt]{article}
\usepackage[T2A]{fontenc}
\usepackage[cp1251]{inputenc}
\begin{document}

\author{S.V. Ludkovsky.}

\title{CW-groups associated with wrap groups.}

\date{12 May 2008}
\maketitle

\begin{abstract}
This article is devoted to the investigation of wrap groups of
connected fiber bundles. CW-groups associated with wrap groups are
studied.
\end{abstract}

\section{Introduction.}
\par Wrap groups of quaternion and octonion fibers
as well as for wider classes of fibers over $\bf R$ or $\bf C$ were
defined and various examples were given together with basic theorems
in \cite{luwrgfbqo}. Studies of their structure were begun in
\cite{lulaswgof,lulatcwg}. This paper continues previous works of
the author on this theme, where generalized loop groups of manifolds
over $\bf R$, $\bf C$ and $\bf H$ were investigated
\cite{ludan,lugmlg,lujmslg,lufoclg}.
\par In this article a structure of wrap groups as CW-groups is
studied. Here the notations and definitions and results from earlier
papers
\cite{luwrgfbqo,lulaswgof,lulatcwg,ludan,lugmlg,lujmslg,lufoclg} are
used.

\section{CW-groups for wrap groups}
\par To avoid misunderstandings we first give our definitions and
notations.

\par {\bf 1. Definitions.} Suppose that $K$ is a Hausdorff space,
which is a union of disjoint open cells, denoted by $\sf e$, ${\sf
e}^n$, ${\sf e}^n_j$, satisfying the following conditions.
\par The closure ${\bar {\sf e}}^n$ of each $n$-cell, ${\sf e}^n\in
K$, is an image of $n$-simplex $\sigma ^n$, in a mapping $f: \sigma
^n\to {\bar {\sf e}}^n$ such that
\par $(CW1)$ $f|_{\sigma ^n\setminus \partial \sigma ^n}$ is a
homeomorphism onto ${\sf e}^n$;
\par $(CW2)$ $\partial {\sf e}^n\subset K^{n-1}$, where $\partial
{\sf e}^n=f(\partial \sigma ^n)={\bar {\sf e}}^n\setminus {\sf
e}^n$, $K^{n-1}$ is the $(n-1)$-dimensional section of $K$
consisting of all cells whose dimensions do not exceed $(n-1)$, in
another words a $(n-1)$-skeleton, $K^{-1} := \emptyset $. Then $K$
is called a cell complex or a complex.
\par Such mapping $f: \sigma ^n\to {\bar {\sf e}}^n$ is called
a characteristic mapping for ${\sf e}^n$.
\par A sub-complex $L\subset K$ is the union of a subset of cells
of $K$, which are cells of $L$, so that if ${\sf e}\subset L$, then
${\bar {\sf e}}\subset L$. If $X$ is a subset of points in $K$, then
$K(X)$ denotes the intersection of all sub-complexes of $K$
containing $X$.
\par A complex $K$ is called closure finite if and only if
$K({\sf e})$ is a finite sub-complex for each cell ${\sf e}\in K$.
\par A weak topology in $K$ is characterized by the condition:
a subset $X$ is closed (or open) in $K$ if and only if $X\cap {\bar
{\sf e}}$ is closed (or relatively open correspondingly) for each
cell ${\sf e}$ of $K$.
\par By a CW-complex we mean one which is closure finite and has the
weak topology.
\par A mapping $f: K\to L$ for CW-complexes $K$ and $L$ is called
cellular, if $f(K^n)\subset L^n$ for each $n=0,1,2,...$.
\par A topological group is called a CW-group if it is a CW-complex
such that the inversion and product mappings $G\ni g\mapsto
g^{-1}\in G$ and $G\times G\ni (g,f)\mapsto fg\in G$ are both
cellular, that is, they carry the $k$-skeleton into the
$k$-skeleton. Then a CW-group $G$ is called countable, if it is a
countable CW-complex.
\par A mapping $f: X\to Y$ is called a homotopy equivalence, if
and only if it has a homotopy inverse meaning a mapping $g: Y\to X$
such that $gf\approx 1_X$ and $fg\approx 1_Y$ (see
\cite{swit,whiteh}).
\par Denote by $({\cal P}^ME;y_0,y_1)_{t,H}$ the quotient
uniform space of $R_{t,H}$ equivalence classes of $H^t_p$ mappings
of a parallel transport structure ${\bf P}_{{\hat {\gamma }},u}$
from $\hat M$ into $E$ such that ${\hat \gamma }: {\hat M}\to N$,
$E=E(N,G,\pi ,\Psi )$ is a principal fiber bundle with a structure
group $G$, $\Xi : {\hat M}\to M$ is a quotient mapping, ${\hat
{\gamma }} ({\hat s}_{0,q})=y_0$, ${\hat {\gamma }}({\hat
s}_{0,q+k})=y_1$ for each $q=1,...,k$. Recall that the equivalence
relation $R_{t,H}$ is generated by: $f\sim g$ if and only if there
exists sequences $f_n$ and $g_n$ converging to $f$ and $g$
respectively in $H^t_p({\hat M},W)$ when $n$ tends to the infinity
such that $f_n=g_n\circ \psi _n$, $\psi _n$ is an
$H^t_p$-diffeomorphism of $\hat M$ preserving marked points ${\hat
s}_{0,j}$, $j=1,...,2k$ (see \S \S 1-3 \cite{luwrgfbqo}). \par We
call $({\cal P}^ME;y_0,y_1)_{t,H}$ the quotient path space.
Particularly, may be $G=e$, that is $E=N$ is a manifold for $G=e$.
As usually consider arcwise connected $E$, $N$ and $G$, where $G$ is
a Lie either alternative or associative group.

\par {\bf 2. Theorem.} {\it If $N$ and $\hat M$ are compact connected
Riemannian $C^{\infty }$ manifolds may be with corners such that the
Ricci tensor $R_{k,l}$ of $N$ is everywhere positive definite, then
the quotient path space $({\cal P}^MN;y_0,y_1)_{t,H}$ for marked
points $y_0$ and $y_1$ in $N$ has the homotopy type of a CW-complex
having only finitely many cells in each dimension.}

\par {\bf Proof.} Theorem A in \cite{milmorse} states if $X$ is
the homotopy direct limit of $ \{ X_j \} $ and $Y$ is the homotopy
direct limit of $ \{ Y_j \} $, if also $f: X\to Y$ is a continuous
map that carries each $X_j$ into $Y_j$ by a homotopy equivalence,
then $f$ itself is a homotopy equivalence. The corollary on page 153
from Theorem A \cite{milmorse} states that if $X$ is the homotopy
direct limit of $ \{ X_j \} $ and each $X_j$ has the homotopy type
of a CW-complex, then $X$ itself has the homotopy type of a
CW-complex. In particular, the quotient space relative to a
continuous quotient mapping of a CW-complex has the homotopy type of
a CW-complex. Therefore, it is sufficient to prove this theorem for
the path space $(P^{\hat M}N;y_0,y_1)_{t,H} := \{ f\in H^t_p({\hat
M},W): \pi \circ f({\hat s}_{0,q})=y_0, \pi \circ f({\hat
s}_{0,q+k})=y_1 ~ \forall q=1,...,k \} $.

\par Since $t\ge [dim (M)/2] + 1$, while $\hat M$ and $N$ are $C^{\infty }$
manifolds, then $C^0\subset H^t_p$ due to the Sobolev embedding
theorem and the homotopy type of $(P^{\hat M}N;y_0,y_1)_{t,H}$ is
the same as $(P^{\hat M}N;y_0,y_1)_{\infty ,H}$.
\par The manifold $N$ is compact, hence it is finite dimensional and the
space consisting of all vectors $v$ of the unit length on $N$ is
compact. The Ricci tensor is the bilinear pairing $R: T_yN\times
T_yN\to \bf R$, which is the trace of the linear transformation
$w\to {\hat R}(v_1,w)u_2$ from $T_yN$ into $T_yN$, where $\hat R$
denotes the Riemann curvature tensor and $R$ is its contraction.
Therefore, there exists $\min \{ R(v,v): v\in T_yN, y\in N, \| v \|
=1 \} =: (n-1)\rho ^{-2}$, where $n$ denotes the dimension of $N$.
\par The manifold $\hat M$ is compact, consequently, there exists a
finite partition $\cal T$ of $\hat M$ consisting of $U_j$ such that
each $U_j$ is homeomorphic with a cube $[0,1]^m$, while
$U_j\setminus
\partial U_j$ is $C^{\infty }$ diffeomorphic with $[0,1]^m\setminus
\partial [0,1]^m$, $\bigcup_j U_j = \hat M$, $m$ denotes the dimension of
$\hat M$, $U_j\cap U_l=
\partial U_j\cap \partial U_l$, $j=1,...,a_0$, $a_0\in \bf N$.
\par Consider a path
${\hat \gamma }: {\hat M}\to N$ such that ${\hat \gamma }({\hat
s}_{0,q})=y_0$ and ${\hat \gamma }({\hat s}_{0,q+k})=y_1$ for each
$q=1,...,k$, where $\hat M$ is the corresponding $C^{\infty }$
Riemannian manifold satisfying Conditions \S 2 \cite{luwrgfbqo} and
$\Xi : {\hat M}\to M$ is the quotient mapping as in \S 2
\cite{luwrgfbqo}, $\Xi ({\hat s}_{0,q})= \Xi ({\hat
s}_{0,q+k})=s_{0,q}$ for each $q=1,...,k$, $s_{0,q}$ and ${\hat
s}_{0,q}, ~ {\hat s}_{0,q+k}$ are marked points in $M$ and $\hat M$
respectively for every $q=1,...,k$, $k\in \bf N$. Therefore, the
path $\hat \gamma $ can be presented as the combination of its
restrictions ${\hat \gamma }|_{U_j}$. \par Without loss of
generality we can take a partition $\cal T$ such that each marked
point ${\hat s}_{0,q}$ in $\hat M$ belongs to $\bigcup_{j=1}^{a_0}
\partial U_j$. If $U_j$ has less, than two distinct marked points
$s_{0,q}$, then introduce in $U_j$ additional marked points
$x_{0,a,j}$ such that to have not less than two distinct marked
points in $U_j$. The manifold $N$ has the homotopy type of a
CW-complex, hence $N^b$ has the homotopy type of a CW-complex for
each $b\in \bf N$ (see also \cite{beggs,milcw} and below).
\par In view of the Sard theorem II.2.10.2 \cite{dubnovfom} and
\S III.6 \cite{miha} the set of all $H^t_p$ diffeomorphisms of $\hat
M$ is everywhere dense in the uniform space $H^t_p({\hat M},{\hat
M})$. Then $(P^{\hat M}N;y_0,y_1)_{t,H}$ has the homotopy type of
$(\bigcup_{j=1}^{a_0} (P^{U_j}N;y_{0,j},y_{1,j})_{t,H})\times
N^{2a_0-2}$, where $y_{0,j}, ~ y_{1,j}$ are $2a_0$ distinct marked
points in $N$ containing $y_0, y_1$ with the corresponding marked
points in $U_j$.

\par In accordance with Proposition $(H)$ \cite{whiteh} if $L$ is a
locally finite complex and $K$ is a CW-complex, then $K\times L$ is
a CW-complex.
\par The sum of CW-complexes is a
CW-complex, the product of CW-complexes is a CW-complex in
accordance with Section 5 and Proposition $(H)$ of \cite{whiteh}.
The manifolds $\hat M$ and $N$ are connected, consequently, it is
sufficient to prove this theorem in the special case of ${\hat M}
=[0,1]^m$.

\par Therefore, consider ${\hat \gamma }: [0,1]^m\to N$, ${\hat \gamma }
(x)\in N$,
$x=(x_1,...,x_m)$, $x_j\in [0,1]$ for each $j=1,...,m$. Suppose that
$\eta _s(x_s)$ is a geodesic between points $a_s$ and $b_s\in N$,
where $\eta _s(x_s) := \eta (z_1,...,z_{s-1},x_s,z_{s+1},...,z_m)$
with marked values of $z_1,...,z_{s-1},z_{s+1},...,z_m\in [0,1]$ and
$\eta : [0,1]^m\to N$, $a_s=\eta _s(0)$, $b_s=\eta _s(1)$. If $\eta
_s(x_s)$ has a length greater than $\pi \rho $, then it has an index
$\lambda \ge 1$ (see also \S \S 16, 17, 19 in \cite{milmorse}).
\par Let ${\bf E}(\zeta )$ denotes the energy functional
of a geodesic in the Riemannian manifold and ${\bf E}_{**}$ be its
Hessian (see \S 12 in \cite{milmorse}).
\par Generally consider a geodesic $\zeta $ of length greater than
$g\pi \rho $, consequently, $\zeta $ has an index $\lambda \ge g$,
where $g\in \bf N$. For each $j=1,...,g$ there exists a vector field
$Y_j$ in $N$ such that $Y_j$ along $\zeta $ vanishes outside the
interval $((j-1)/k, j/k)$, and so that ${\bf E}_{**}(Y_j,Y_l)<0$.
Since ${\bf E}_{**}(Y_j,Y_l)=0$ for each $j\ne l$, then $Y_1,...Y_g$
span a $g$-dimensional subspace of $\bigcup_{y\in \zeta ([0,1])}
T_yN$ on which ${\bf E}_{**}$ is negative definite (see \S 19 in
\cite{milmorse}).

\par Suppose that points $y_{0,j}$ and $y_{1,j}$ are not conjugate
along any geodesic from $y_{0,j}$ to $y_{1,j}$, hence there exists
only a finite number of geodesics like $\eta _s$ from $y_{0,j}$ to
$y_{1,j}$ in $N$ by the variable $x_s$ of length not greater than
$g\pi \rho $. Hence there exists only finitely many geodesics with
index less than $g$.
\par In accordance with Theorem 17.3 \cite{milmorse}
if $N$ is a complete Riemannian manifold and $y_0, y_1\in N$ are two
points, which are not conjugate along any geodesic, then
$(P^{[0,1]}N;y_0,y_1)_{t,H}$ has the homotopy type of a countable
CW-complex containing one cell of dimension $\lambda $ for each
geodesic from $y_0$ to $y_1$ of index $\lambda $.
\par Together with Theorem 17.3 \cite{milmorse} this completes
the proof for $dim (M)=1$. For $m>1$ proceed by induction:
\par $(P^{[0,1]^m}N;y_0,y_1)_{\infty ,H} = (P^{[0,1]^{m-1}}
(P^{[0,1]}N;y_0,y_1)_{\infty ,H};y_0,y_1)_{\infty ,H}$, where $y_b$
in $(P^{[0,1]^l}N;y_0,y_1)_{\infty ,H}$ denotes the constant mapping
$y_b: [0,1]^l\to N$, $y_b([0,1]^l) = \{ y_b \} $, $ \{ y_b \} $
denotes the singleton in $N$, $b=1, 2$, $l\in \bf N$, here the
notation $y_b$ corresponds to $y_{b,j}$ for some $j$.
\par This procedure lowers a number of variables on each step
by one.  In view of Theorem 19.6 \cite{milmorse}
$(P^{[0,1]}N)_{t,H}$ has the homotopy type of a CW complex $B$,
which is $\sigma $-compact, that is a countable union of compact
sets.
\par Consider now $(P^{[0,1]}B)_{t,H}$, where $B$ is a countable union
of compact Riemannian manifolds may be with corners, since each
polyhedron in $\bf R^n$ with $n\in \bf N$ is a manifold with
corners. Put $B = \bigcup_{j\in \Lambda }B_j$, $B^k :=
\bigcup_{j=1}^k B_j$, where $B_j$ is a compact Riemannian manifold
with corners being a $j$-skeleton of a CW complex, $\Lambda \subset
\bf N$. Up to a homotopy type or bending $B_j$ a little in the
corresponding Euclidean space $\bf R^n$ of dimension $n\ge 2~ dim
~(B_j)$, $B_j\hookrightarrow {\bf R}^j\hookrightarrow {\bf R}^n$, we
can consider, that each $B_j$ is homotopy equivalent to a compact
Riemannian manifold $X_j$ with positive definite Ricci tensor.
Therefore, we have to consider now $(P^{[0,1]}X)_{t,H}$, where
$X=\bigcup_j X_j$. Put $X^j=\bigcup_{k\le j} X_k$, then $X^j\subset
X^{j+1}$ for each $j\in \Lambda $, $dim ~(X_j)=j$.
\par Each path from the
compact manifold $M$ into a CW-complex $B$ has a compact image,
consequently, it has a finite covering by cells. Hence a continuous
path from $M$ into $X$ up to a homotopy equivalence has a finite
covering by $X^j$.
\par If $N_1$ and $N_2$ are homotopy equivalent Riemannian
manifolds, then \\ $(P^{[0,1]}N_1;y_{0,1},y_{1,1})_{t,H}$ and
$(P^{[0,1]}N_2;y_{0,2},y_{1,2})_{t,H}$ are homotopy equivalent, when
$y_{0,1}\ne y_{0,2}$ and $y_{0,2}\ne y_{1,2}$ simultaneously. On the
other hand, $(P^{[0,1]}X;y_0,y_1)_{t,H}$ is homotopy equivalent with
a CW-complex $K = \bigcup_{j\in \Lambda } K_j$, where each $K_j$ is
a CW-complex homotopy equivalent with
$(P^{[0,1]}X^j;y_0,y_1)_{t,H}$, where $y_0, y_1\in X_1$, so that
$K_j\subset K_{j+1}$ for each $j$, since $X^j\subset X^{j+1}$.
\par Denote by $\cal W$ the class of all spaces having
the homotopy type of a CW-complex. By a CW-n-ad ${\sf K} =
(K;K_1,...,K_{n-1})$ is undermined a CW-complex together with
$(n-1)$ numbered sub-complexes $K_1,...,K_{n-1}$. Then ${\cal W}^n$
denotes the class of all $n$-ads which have the homotopy type of a
CW-$n$-ad. As usually ${\sf A}^{\sf C}$ denotes the subspace of the
space $A^C$ of all continuous functions $f$ from $A$ into $C$ such
that $f: {\sf C}\to {\sf A}$ is a mapping of $n$-ads, that is the
induced mappings are $f_j: C_j\to A_j$ from the $j$-skeleton to the
$j$-skeleton for each $1\le j\le n$.

\par In accordance with Theorem 3 \cite{milcw} if $A$ belongs to
the class ${\cal W}^n$ and $\sf C$ is a compact $n$-ad, then the
function space ${\sf A}^{\sf C}$ belongs to $\cal W$. In fact the
$n$-ad $(A^C;(A,A_1)^{(C,C_1)},...,(A,A_{n-1})^{(C,C_{n-1})})$
belongs to the class ${\cal W}^n$.
\par Thus, $(P^{\hat M}N;y_0,y_1)_{t,H}$ has the homotopy type
of the CW-complex.

\par {\bf 2. Corollary.} {\it If $M$ and $\hat M$ and $N$ are manifolds
$H^t_p$ and $H^{t'}_p$ diffeomorphic with $C^{\infty }$ Riemannian
manifolds $M_1$ and ${\hat M}_1$ and $N$ correspondingly, $t'\ge t$,
where $M_1$, ${\hat M}_1$ and $N_1$ satisfy conditions of the
preceding theorem, then the path space $(P^{\hat M}N;y_0,y_1)_{t,H}$
and the quotient path space $({\cal P}^MN;y_0,y_1)_{t,H}$ for marked
points $y_0$ and $y_1$ in $N$ are of the homotopy types of
CW-complexes having only finitely many cells in each dimension.}
\par {\bf Proof.} Let $\phi : {\hat M}_1\to \hat M$ and
$\theta : N_1\to N$ be homeomorphisms, which are $H^t_p$ and
$H^{t'}_p$ diffeomorphisms. Then the uniform spaces $(P^{\hat
M}N;y_0,y_1)_{t,H}$ and $(P^{{\hat M}_1}N_1;y_{0,1},y_{1,1})_{t,H}$
are isomorphic, where the mapping $f\mapsto \theta ^{-1}\circ f\circ
\phi $ establishes the isomorphism, $f\in (P^{\hat
M}N;y_0,y_1)_{t,H}$, $\theta (y_{b,1})= y_b$ for $b=1, 2$. Using
this isomorphism and applying the preceding theorem to $(P^{{\hat
M}_1}N_1;y_{0,1},y_{1,1})_{t,H}$ and the quotient path space $({\cal
P}^{M_1}N_1;y_{0,1},y_{1,1})_{t,H}$ we get the statement of this
corollary.

\par {\bf  3. Corollary.} {\it Let $M$ and $N$ be satisfying
conditions of the preceding Corollary. Then the wrap monoid
$(S^MN)_{t,H}$ and the wrap group $(W^MN)_{t,H}$ have homotopy types
of CW-complexes having only finitely many cells in each dimension.}
\par {\bf Proof.} The wrap monoid has the homotopy type of
$({\cal P}^MN;y_0,y_0)_{t,H}$. On the other hand, the wrap group is
the quotient of the free commutative group $F$ generated by
$(S^MN)_{t,H}$ by the closed equivalence relation, which is obtained
factorizing by the minimal closed normal subgroup $B$ containing all
elements of the form $[a+b]-[a]-[b]$, where $a, b\in (S^MN)_{t,H}$,
$[a]$ and $[b]$ are the corresponding elements of $F$. Topologically
$F$ is isomorphic with $[(S^MN)_{t,H}]^{\bf Z}$ supplied with the
weak (Tychonoff) product topology. Applying Corollary on page 153
from Theorem A \cite{milmorse} and the preceding theorem we get the
statement of this corollary.

\par {\bf 4. Corollary.} {\it Let $M$ and $N$ be satisfying
conditions of Corollary 2, while $E$ be a principal fibre bundle
with the structure group $G$, which is up to the homotopy a
CW-group. Then a wrap monoid $(S^ME)_{t,H}$ and a wrap group
$(W^ME)_{t,H}$ have homotopy types of a CW-monoid and a CW-group
correspondingly.}
\par {\bf Proof.} By Proposition $(N)$ any covering complex of a
CW-complex is a CW-complex \cite{whiteh}. Therefore, if prove that
$(S^ME)_{t,H}$ is a CW-complex, then it would mean that
$(W^ME)_{t,H}$ is a CW-complex. This follows immediately from the
preceding corollary and Proposition 7.1 \cite{lulaswgof} and
Proposition $(H)$ \cite{whiteh}, since $(S^ME)_{t,H}$ and
$(W^ME)_{t,H}$ have structures of principal $G^k$-bundles over
$(S^MN)_{t,H}$ and $(W^MN)_{t,H}$.
\par On the other, hand the mapping $(S^MN)_{t,H}\ni (f,g)\to fg\in
(S^MN)_{t,H}$ is cellular, since if $a, b\in K^n$, then $a\vee b\in
K^n\vee K^n$, where the bunch $K^n\vee K^n$ of $K^n$ by a finite
number of marked points consists of cells of dimension at most $n$.
Therefore, in $(W^MN)_{t,H}$ the group multiplication is cellular as
well (see also \S 3). In $(W^MN)_{t,H}$ the mapping $f\mapsto
f^{-1}$ is cellular due to the definition of the wrap group. Since
$G$ is the CW-group, then $G^k$ is the CW-group, consequently,
$(S^ME)_{t,H}$ and $(W^ME)_{t,H}$ are the CW-monoid and the CW-group
respectively.

\par {\bf 5. Remark.} A topological space $P$ is said to be dominating
a topological space $X$ if and only if there are continuous mappings
$f: X\to P$ and $g: P\to X$ such that $gf\approx 1_X$. In accordance
with Theorem 1 \cite{milcw} $A$ belongs to the class ${\cal W}_0$ if
and only if $A$ is dominated by a countable CW-complex.
\par  If $G$ is a compact simply connected Lie
group, then in accordance with Theorem 21.7 \cite{milmorse}
$(P^{[0,1]}G;y_0,y_1)_{t,H}$ has the homotopy type of a CW-complex
with no odd-dimensional cells and with only finite number of
$n$-cells for each even number $n$. These two theorems imply that
$G$ also is a CW-group, since $(P^{[0,1]}G)_{t,H}$ dominates $G$ and
applying the homotopy equivalence.
\par If $G$ is not associative, but alternative, then the
corresponding CW-group is alternative as well, since if $a_1\approx
a_2$, $b_1\approx b_2$ are homotopy equivalent elements of $G$, then
$(a_1a_1)b_1 = a_1(a_1b_1) \approx a_1(a_2b_2) \approx a_2(a_2b_2) =
(a_2a_2)b_2$ and $b_1=a_1^{-1}(a_1b_1) \approx a_1^{-1}(a_2b_2)
\approx a_2^{-1}(a_2b_2) = (a_2^{-1}a_2)b_2=b_2$ and analogously for
identities with $a_j$ on the right from $b_j$.
\par In accordance with Corollary 1
\cite{milcw} every separable finite dimensional manifold belongs to
the class ${\cal W}_0$, where ${\cal W}_0$ denotes the class of
topological spaces having the homotopy type of countable
CW-complexes. Due to Corollary 2 \cite{milcw} if $A$ belongs to
${\cal W}_0$ and $C$ is a compact metric space, then the function
space $A^C$ in the compact open topology belongs to ${\cal W}_0$.
Therefore, modifying Theorem 2 and Corollary 4 we get.
\par {\bf 6. Proposition.} {\it If $N$ is a finite dimensional
separable manifold, $G$ is a CW-group, then $({\cal
P}^ME;y_0,y_1)_{t,H}$ has the homotopy type of a CW-complex,
$(S^ME)_{t,H}$ and $(W^ME)_{t,H}$ have homotopy types of a CW-monoid
and a CW-group respectively.}

\end{document}